\documentclass{amsart}
\usepackage{times,amsfonts,mathabx,amssymb,amsrefs,mathrsfs,tikz}
\usepackage[printwatermark]{xwatermark}
\usepackage{xcolor}
\usepackage{graphicx}

\usetikzlibrary{decorations,decorations.pathmorphing}

\newcommand{\Z}{{\mathbb Z}}

\newcommand{\constanta}{{\rm Const}}

\newtheorem{thm}{Theorem}

\newtheorem{cor}[thm]{Corollary}
\newtheorem{lem}{Lemma}[section]
\theoremstyle{definition}
\newtheorem{rem}{Remark}[section]
\newtheorem{convention}{Convention}[section]

\newtheorem{corlem}{Corollary}[section]
\newtheorem{exa}[thm]{Example}

\title[Iterated identities]{No iterated identities satisfied by all finite groups}

\author{Alexei Belov}
\address{ 
  Departement of mathematics, Bar Ilan University, Israel}
\email{kanelster@gmail.com }
\author{Anna Erschler}
\address{D\'{e}partement de math\'{e}matiques et applications, \'{E}cole normale
sup\'{e}rieure, CNRS, PSL Research University, 45 rue d'Ulm, 75005 Paris}
\email{anna.erschler@ens.fr}

\date{\today} 
\thanks{The work of the authors is partially supported by the ERC grant GroIsRan; the work of the first named author is also supported the
Israel science foundation Grant 1623/16.}

\begin{document}

\maketitle

\begin{abstract}
We show that there is no iterated identity satisfied by all  finite groups.
For $w$ being a non-trivial word of length $l$, we show that there exists a finite group $G$ of cardinality at most $\exp(l^C)$ which does not satisfy the iterated identity $w$.
The proof uses the approach of Borisov and Sapir, who used dynamics of polynomial mappings for the proof of non residual finiteness of some groups.

\end{abstract}

\section{Introduction}

It is well-known and not difficult to see that there is no non-trivial group identity which is satisfied by all finite groups. We strengthen this fact by showing that there is no {\it iterated} group identity which is satisfied by all finite groups, and we construct a group violating a given iterated identity, providing an upper bound for the cardinality of this group.
We recall the definition of iterated identity from \cite{erschleriterated}. 
We say that a group $G$ satisfies {\it an Engel type iterated identity} $w$  if
for any $x_1, \dots, x_m \in G$ there exists $n$ such that 
\begin{equation} \label{eq:iterated}
  w_{\circ n}(x_1, x_2, \dots, x_m)  =w(w(\dots(w(x_1, x_2, \dots, x_m), x_2, \dots, x_m), x_2, \dots, x_m))=e.
\end{equation}
In the sequel, we call Engel type iterated identities for short {\it iterated identities}. For definitions other than that of Engel type see \cite{erschleriterated}.

The fact the group satisfies an iterated identity depends only on the element of the free group represented by this word, in other words, the property to satisfy an iterated identity does not change if we replace a word by a freely equivalent one, in particular, any group satisfies $w$ if $w$ is freely equivalent to an empty word.

The definition of iterated identities is close to the notion of "correct sequences", studied by Plotkin, Bandman, Greuel, Grunewald, Kuniavskii, Pfister, Guralnick and Shalev  in 
\cite{plotkin,bandmanetal,guralnickplotkinshalev}. Examples of such sequences, without this terminology, were previously constructed by Brandl and Wilson \cite{brandlwilson}, Bray and Wilson \cite{braywilson}
and Ribnere to characterize finite solvable groups. See also \cite{grunewaldkuplotkin}.

For some groups and some classes of groups a priori not bounded number of iteration in the definition of iterated identity is essential, as for example it is the case for the first Grigorchuk group, which is  a $2$ torsion group \cite{grigorchuk}, that is, it satisfies the iterated identity $w(x_1)=x_1^2$, but this group does not satisfy any identity by a result of Abert  \cite{abert}. For some other groups the number of such iterations is bounded for all $x_1$, \dots, $x_n$, a strong version of this phenomena is when such bound does not depend on the iterated identity $w$, as it is for example the case for 
any finitely generated metabelian group \cite{erschleriterated}.

\begin{thm} \label{thm:noid}
Let $w(x_1, \dots, x_m)$ be a word (on $n$ letters, $m\ge1$) which is not freely equivalent to an empty word.
Then there exists a finite group $G$ such that $G$ does not satisfy an iterated identity $w$. 

Moreover, there exists $C>0$ such that for any $n\ge 1$ and any word  $w$  on $n$ letters the group $G$ can be chosen to have at most  $ \exp(l^{C})$ elements, where $l=l(w)$ is the length of the word $w$. 
\end{thm}

An upper bound for the cardinality of a finite group in the second part of the theorem might not be optimal. One can ask whether one can replace $\exp(l^{C})$ by $l^C$. For related questions see also Section
\ref{se:openquestions}.

A standard argument to show that there is no identity for all finite groups is to observe that free non-Abelian groups are  residually finite, and to conclude that if $w$ is an identity satisfied by all finite groups, then the free group $F_2$ also satisfies $w$. Observe that this argument does not work for iterated identities.
Indeed, free groups are residually nilpotent, however every nilpotent group satisfies the iterated identity $w(x_1, x_2)= [x_1, x_2]$, while a free group does not satisfy any non-trivial iterated identity.

To prove the theorem, we show that for any word $w$ on $x_1$, $x_2$ representing an element in the commutator subgroup of $F_2$ there exists $n \ge 1$ such that 
$w_{\circ n}(x_1, x_2) =x_1$ admits a solution with $x_1\ne 1$ in some finite group.  Here $w_{\circ n}$ is as defined in the equation $1$.

 Much progress has been achieved in  recent years in understanding the image of the verbal mapping from $G^n \to G$
$(x_1, \dots, x_n) \to w(x_1, \dots, x_n)$. Larsen, Shalev and Tiep prove in \cite{larsenshalevtiep} that for any word $w$ and for any sufficiently large finite simple non-Abelian group $w(G^n)w(G^n)= G$, that is, for any $g\in G$ there exists $x_1, \dots, x_n \in G$ and $x'_1, \dots, x'_n \in G$ such that $w(x_1, \dots, x_n) w(x'_1, \dots, x'_n) =g$. Moreover, for some words $w$ such verbal mapping turn out to be surjective. 
Libeck,  O'Brien , Shalev and Tiep  \cite{liebecketal}, proving  the Ore conjecture, show that this is the case for $w(x_1,x_2)=[x_1,x_2 ]$ and any finite simple non-Abelian group. Observe however that the image of the mappings from $G^n \to G^n$ which sends 
$(x_1, \dots, x_n)$ to  $\left( w(x_1, \dots, x_n), x_2, x_3, \dots, x_n \right)$ is far from being surjective, and the structure of periodic points for such mappings seems to be less understood.

To solve the equation $w_{\circ m}(x_1, x_2) =x_1$, we use the idea and the result of Borisov and Sapir from \cite{borisovsapir}, who use {\it quasi-periodic} points of polynomial mappings to prove non-residual finiteness of some one relator groups, namely
of what is called {\it mapping tori} (also called ascending $HNN$ extensions) of injective group endomorphisms: those are groups of the form $(x_1, x_2, \dots , x_k,  t |   R,  tx_it^{-1}=w_i, i \le i \le k)$, where $x_i \to w_i$ is an injective endomorphism of
the group $(x_1, x_2, \dots, x_i  |  R)$.

In contrast with mapping tori of groups endomorphisms, general one relator groups are not necessary residually finite, and it is a long
standing problem to characterize residually finite one relator groups. A conjecture of Baumslag, proven by Wise in \cite{wise} states that  one relator  group containing a non-trivial torsion
element is residually finite. The situation for groups without torsion elements is less understood.

Consider a sequence of the one relator group $G_m=[x_1,x_2: w_{\circ m}(x_1, x_2) =x_1]$. If a finite quotient of a group $G_m$ is such that  the image  of $x_1$ is not equal to one  in this finite quotient, then
in this finite quotient  the image of $x_1$ is a non-fixed periodic point for the verbal map $x \to w(x,y)$, for a fixed $y$. 

We will construct finite quotients of groups $G_m$  as  subgroups of  $SL(2,\mathcal{K})$, for an appropriately chosen finite field $\mathcal{K}$.  In  Section \ref{se:moduloBS}
 we  outline the proof of the theorem and prove its first claim. To do this, we choose a two-times-two integer valued matrix $y_0$ which can be one of the free generators of a free non-Abelian subgroup in  $SL(2, \mathbb{C})$, regard $w(x,y_0)$ as a function of $x$, observe that the entries of $w(x,y_0)$ are rational functions $R_{i_1, i_2}(x_{1,1}, x_{1,2},x_{2,1},x_{2,2})$ in $x_{1,1}$, $x_{1,2}$, $x_{2,1}$, $x_{2,2}$. Multiplying by a power of the determinant  of the matrix for $x$, we get 
 polynomials $H_{i_1, i_2}(x_{1,1}, x_{1,2},x_{2,1},x_{2,2})$. We will  need to check that the system of 
  equations $H_{i_1, i_2}(x_{1,1},x_{1,2}, x_{2,1}, x_{2,2})= x_{i_1,i_2}^Q$ satisfies the assumption of Theorem $3.2$ of \cite{borisovsapir} and we aplly this theorem  to solve this system of equations,
assuming that $q$ is a large enough prime and  $Q$ is a large power of $q$.
 
We check that the image on the $4$-th iteration of the polynomial mapping in question contains at least  one point with non-zero determinant, and such that the matrix is not a diagonal matrix. Hence we obtain at least one
 non-trivial solution of the system of the equations for $H_{i_1,i_2}$, which is non-diagonal matrix with non-zero determinant. Normalizing, if necessary, this solution by the square of the determinant of the corresponding matrix, we will
obtain a non-identity solution for the system of the equations $R_{i_1, i_2}(x_{1,1},x_{1,2}, x_{2,1}, x_{2,2})= x_{i_1,i_2}^Q$.  This solution belongs to a finite extension of $F_q$. Such solution 
 provides a  non-identity periodic point  $x$ for iteration of $w$, for some $m\ge 1$, and this implies in particular that
$w$ is not an iterated identity in the subgroup generated by $x$ and $y_0$ in $SL(2,k)$.

In Section \ref{se:secondpart} we obtain a bound for the cardinality of a subgroup $SL(2, \mathcal{K})$. 
To do this, we need to control  the cardinality of the finite field $\mathcal{K}$. For this purpose, instead of using Theorem $3.2$ of  of \cite{borisovsapir}, we  prove and use
a
version of that theorem, see Theorem  \ref{thm:borisovsapir}.  Given $n$ polynomials  $f_i$ on $n$ variables over a finite field,  and a polynomial $D_0$,
this theorem provides a lower bound for $Q$ in terms of degree of these polynomials with the following property.
If $D_0$ is equal to zero on any solution over algebraic closure of $F_q$ of the system of equations
$$
f_i(x_1, \dots, x_n)= x_i^Q,
$$
then $D$ is equal to zero on any point of the $n$-th iteration of the polynomial mapping $f=(f_1, \dots, f_n)$.
To prove this theorem  we follow the strategy of the proof of Borisov Sapir, 
the main ingredient of the proof is Lemma 3.4, which  provides an explicit
estimate for  Lemma 3:5 in [11].  Given a solution $a_1, \dots, a_n$ of the system of the equations, this lemma provides an estimate for $k$ such that $(f_i^{(n)}-\constanta)^k$  belongs to the localisation at
$(a_1, a_2, \dots, a_n)$ of the ideal generated by $f_i(x_1, \dots, x_n)- x_i^Q$,  for each $i$. Here $f_1^{(n)}, \dots, f^{(n)}_n$ is the $n$-th iteration of the polynomial mapping $f$.
As the last step of the proof, rather than using one of two possible arguments used in \cite{borisovsapir}, we make use of  the fact that the polynomials $H^{(4)}_i- x_i^Q$ form
a {\it Gr{\"o}bner basis} with respect to Graded Lex order.

In section \ref{se:generaldynamics} we give a more general version of Theorem \ref{thm:noid}, where instead of iterations on one variable we consider iterations of verbal mappings on several variables.
Given an endomorphism  $\phi$  of a free group $F_n$,  $\phi_{\circ m}$ denotes the $m$-th iteration of $\phi$ and  $H_n$  denotes the kernel of $\phi_{\circ n}$. It is clear that $H_n$ is a normal subgroup, and $F_n/H_n$ is isomorphic to the image of $\phi_{\circ n}$, this image is isomorphic to a subgroup of $F_n$, and thus is is finitetly generated free group, $F_n/H_n$.

\begin{thm}  \label{thm:borsapmol} Let $\phi$ be an endomorphism of a free group. 
For any $g \in F_n$, $g\notin H_n$, there exist a finite quotient group $G$ of $F_n$ such that $\phi(g) \ne e$, where $\phi$ is the projection map from $F_n$ to $G$ and such that $\phi$ induces an automorphsim of $G$.

Moreover, we can choose $G$ as above of cardinality at most $\exp(L^{C_n})$, wrere $L= \sum_{i=1}^n \phi(x_i)$, $x_i$ is a free generating set of $F_n$, and $C_n$ is a positive constant depending on $n$.

\end{thm}

\subsection*{Aknowledgements.} The authors are grateful to Boris Kunyavsky, Evgeny Plotkin  and Mark Sapir   for comments on the preliminary
version of the paper.
We would like to thank Andreas Thom for informing us about an impoved estimate
for  the result of \cite{bradfordthom},   Alexander Borisov for
drawing our attention to \cite{varshavsky}.

\section{Idea of the proof of Theorem \ref{thm:noid}  and the proof of its first claim. } \label{se:moduloBS}

We start with a not difficult lemma that shows that it is enough to consider only iterated identities on two letters.
\begin{lem}
Suppose that a class of groups does not satisfy any non-trivial iterated on two letters. Then is no iterated identity satisfied by this class of groups.
\end{lem}

{\bf Proof.} Suppose that $\bar{w}(x_1, x_2, \dots , x_m)$ is an iterated identity and $w$ is not freely equivalent to an empty word.
Choose $u_2 (x,y)$, \dots, $u_m(x,y)$ and put 
$$
w(x_1, x_2) = \bar{w}(x_1, u_2(x_1, x_2),u_m(x_1,  x_2)).
$$
 It is clear that $w$ is an iterated identity.
Now suppose that $u_2$, \dots , $u_n$ are such that $x_1, u_2 (x_1, x_2)$, ...,  $u_m(x_1, x_2)$ generate a free group of rank $m$ in the free group generated by $x_1$ and $x_2$. Observe that in this case $w$ is not freely equivalent an empty word, and thus $w$ is a non-trivial iterated identity on two letters.

Take a non-trivial word $w(x,y)$ on two letters.
The following obvious remark shows that for the proof of the theorem it is enough to consider words representing elements in  the commutator subgroup of $F_2$.

\begin{rem} \label{obviousremark}
If $w$ is a word on $x,y$  which depends on $y$ only, then $w$ is freely equivalent to
$y^m$, $m \ne 1$. If $m\ne 1$ and  $M>m$ is relatively prime with $m$, then $w$ is not iterated identity in a finite cyclic group of $M$ elements.

More generally,  if $w$ is a word on $x,y$  
which does not belong to the commutator group $[F_2,F_2]$, then $w(x,y)=x^my^k\bar{w}$, where at least one of $k$ and $m$ is not equal to $0$,
and $\bar{w}$ is a word representing an element of the commutator group. Considering the iterated values of $w(x,y)$ in $0,y$ and $0,x$ we conclude that
$w$ is not iterated identity in a finite cyclic group of $M$ elements, for any $M$ which is relatively prime with $m$. 
\end{rem}

In the sequel we assume that $w$ is a word on two letters representing an element of the commutator subgroup of $F_2$.

Now consider two times two matrices  $x$ and $y$

\[ x= \left( \begin{array}{ccc}
x_{1,1} & x_{1,2}  \\
x_{2,1} & x_{2,2}
 \end{array} \right) 
\mbox{ and  }
 y= \left( \begin{array}{ccc}
y_{1,1} & y_{1,2}  \\
y_{2,1} & y_{2,2}
 \end{array} \right)\] 

\begin{convention}\label{conv:freesubgroup}

 We assume that  $y_{i_1,i_2} \subset \mathbb{Z}$, $i_1, i_2=1,2$, are such that for some choice  of $x_i$ in $\mathbb{C}$,  the group generated by the matrices $x$ and $y$ is free,
$y$ is in $SL(2, \mathbb{R})$.
For example,  one can take
\[ y= \left( \begin{array}{ccc}
1 & 0  \\
1 & 1
\end{array} \right)\]

\end{convention}

\begin{convention}\label{strongerconv:freesubgroup}
We asume that  $y_i \subset \mathbb{Z}$, $1 \le i \le 4$, are such that for some choice  of $x_i$ in $\mathbb{Z}$,  the group generated by the matrices $x$ and $y$ is free and 
$x,y$ are in $SL(2,\mathbb{Z})$.
For example,  one can take
\[ y= \left( \begin{array}{ccc}
1 & 0  \\
2 & 1
\end{array} \right)\] 

\end{convention}

First observe that  we can chose $y_i$ as in the convention \ref{strongerconv:freesubgroup}, since $SL(2,\mathbb{Z})$ is virtually free, and in particular this group contains free subgroups.

 For example take any integer $m\ge 2$ (e.g. $m=2$), put $\alpha=\beta=m$ and consider

\[ x= \left( \begin{array}{ccc}
1 & \alpha  \\
0 & 1
 \end{array} \right), 
\mbox{ and }
 y= \left( \begin{array}{ccc}
1 & 0  \\
\beta & 1
 \end{array} \right).\]

The subgroup generated by such $x$ and $y$ is free whenever  $\alpha=\beta\ge 2$ are positive integers, see e.g. Theorem $14.2.1$ in \cite{kargapolovmerzljakov}. (For $\alpha=
\beta=2$ this subgroup is called Sanov subgroup).
Moreover, it is easy to see that the subgroup generated by $x$ and $y$ depends up to an isomorphism only on the product $\alpha \beta$ \cite{changetal}; this group
is also free  for any $\alpha$ and $\beta$ such that $\alpha \beta$ is transcendental \cite{fuksrabinovich,changetal} (the group is known to be free for example for any complex $\alpha$, $\beta$ such that
$|\alpha \beta |, |\alpha \beta -2|  >2,  |\alpha \beta +2|$  \cite{changetal}, but apparently it is not known in general when it is free).

In particular, $y$ for $\beta=1$ satisfies the assumption of the Convention \ref{conv:freesubgroup}, since it is sufficient to consider $x$  as above with $\alpha$ which is transcendental.

While for the proof of Theorem 1 any matrix $y$ as in solution $2$ will suffice, in case we want to get more information about finite groups we construct, it might be interesting to start
with various choices of $y$ (for some open questions see Section 5).

Now we fix integers $y_{i_1, i_2}$, $i_1, i_2=1,2$ as in Convention \ref{conv:freesubgroup}. Note that
\[ x^{-1}= \frac{1}{x_{1,1}x_{2,2}-x_{1,2}x_{2,1}}\left( \begin{array}{ccc}
x_{2,2} & -x_{1,2}  \\
-x_{2,1} & x_{1,1}
 \end{array} \right) 
\]

Observe that 
 \[ w(x,y)= \left( \begin{array}{ccc}
R_{1,1}(x_{i_1, i_2},y_{i_1,i_2}) & R_{1,2}(x_{i_1,i_2},y_{i_1,i_2})  \\
R_{2,1} (x_{i_1,i_2},y_{i_1,i_2})& R_{2,2}(x_{i_1, i_2},y_{i_1, i_2}),
 \end{array} \right)\] 
 
\smallskip

where $R_{j_1,j_2}$,  $j_1, j_2 =1,2$ are rational functions in $x_{i_1,i_2},y_{i_1, i_2}$, $i_1, i_2 =1,2$ with integer coefficients. We consider fixed integers $y_{i_1, i_2}$, with $y\in SL(2, \mathbb{Z})$ (e.g. $y_{1,1}=1$, $y_{1,2}=0$, $y_{2,1}=2$, $y_{2,2}=1$) as above, and then
 $R_{j_1,j_2}(x_{i_1,i_2}) = R_{j_1,j_2}(x_{i_1, i_2}, y_{i_1, i_2})$ are rational functions in $x_{1,1}$, $x_{1,2}$, $x_{2,1}$ and $x_{2,2}$ with integer coefficients. For each $j_1, j_2$ it holds 
 
 $$
 R_{j_1, j_2}(x_{1,1}, x_{1,2}, x_{2,1}, x_{2,2}) = H_{j_1, j_2}(x_{1,1}, x_{1,2},
  x_{2,1}, x_{2,2}) /(x_{1,1}x_{2,2}-x_{1,2}x_{2,1})^s,
 $$  where
$H_{j_1, j_2}$ are polynomials with integer coefficients in $x_{1,1}, x_{1,2}, x_{2,1}, x_{2,2}$ and $s$ is the number of occurrences  of $x^{-1}$ in $w$.

Observe that 
$$
H= \left( \begin{array}{ccc}
H_{1,1}(x_{1,1},x_{1,2},x_{2,1}, x_{2,2}) & H_{1,2}(x_{1,1},x_{1,2},x_{2,1}, x_{2,2})  \\
H_{2,1}(   x_{1,1},x_{1,2},x_{2,1}, x_{2,2} )        & H_{2,2}( x_{1,1},x_{1,2},x_{2,1}, x_{2,2}  )
 \end{array} \right) 
$$
\smallskip
is not an identity matrix.

Now observe that if $y$ satisfies the assumption of Convention \ref{strongerconv:freesubgroup},
 then we know moreover that for some values of $x_i \in \mathbb{Z}$ the corresponding matrices $w(x,y)$ and $w(x',y)$
do not commute.  Indeed, observe that if $w(x,y)$ is a freely reduced word on two letters that has at least one entry of $x$ or $x^{-1}$, then $w(x,y)$ and $w(x^k,y)$ do not commute in the free group generated by $x$ and $y$; this implies in particular  that at least one of the rational functions $R_{1,2}$ and $R_{2,1}$ is not zero, and therefore that at least one of polynomials $H_{1,2}$ and $H_{2,1}$ is not zero.

 We consider $Q$ to be a power of $q$ and we want to solve over field $\mathcal{K}$  of characteristic $q$ the system of four equations:

\begin{equation} \label{RQ}
 R_{j_1,j_2}(x_{i_1,i_2}, y_{i_1,i_2}) = x_{j_1,j_2}^Q,
\end{equation}
for $j_1,j_2=1,2$.
To do this, we start by solving the system of polynomials equations:
 \begin{equation} \label{HQ}
 H_{j_1,j_2}(x_{i_1, i_2}, y_{i_1, i_2}) = x_{j_1, j_2}^Q,
\end{equation}
 $j_1, j_2=1,2$.

It is easier to work with the system of the equations (\ref{HQ}) rather then (\ref{RQ}) is that  polynomials $H_{j_1,j_2}(x_{i_1, i_2}, y_{i_1, i_2}) - x_{j_1, j_2}^Q$ form a Groebner basis (in the next section we recall a definition and basic properties of Groebner bases), while polynomials obtained from rational functions
$(R_{j_1, j_2}(x_{i_1, i_2}, y_{i_1, i_2}) = x_{j_1, j_2}^Q)$, after multiplication on the denominator, do not in general form such basis.

The solutions of the system of polynomial equations are Zariski dense in the image the fourth iteration of the polynomial mappings from $\bar{F_q}$ to $\bar{F_q}^4$, where $\bar{F_q }$ denotes the algebraic closure of $F_q$ by  Theorem $3.2$ in Borisov Sapir \cite{borisovsapir}; for a more general statement see
  Corollary $1.2$ on page $5$ of the preprint of  Hrushovski \cite{hrushovski}, see also \cite{vashavsky}.
Indeed, observe we know that the dimension of the fourth iteration of $H$ is not zero, since the image contains at least two points over the field on $q$ elements, for any sufficiently large $q$. (And there exists a variety of dimension greater than $0$ , such that the iteration of the polynomial mapping corresponding to $H$, restricted to this variety,  is dominant).
Moreover, observe that for sufficiently large $q$ there is at least one point $v_{1,1}, v_{1,2}, v_{2,1}, v_{2,2}$  in the image of $f$ such that $v_{1,1}v_{2,2}-v_{1,2}v_{2,1} \ne e$ and either $v_{2,1}$ or $v_{1,2}$ is not equal to $0$. 
Indeed, suppose that $w$ is reduced word containing at least one entry of $x$ or $x^{-1}$.  Take any $x$, $y$  as in Convention \ref{strongerconv:freesubgroup}, that is $x$ and $y$ are  in $SL(2, \mathbb{Z})$ such that $x$ and $y$ generate a free subgroup. Observe that $w(x^m,y)$ belongs to $SL(2,\mathbb{Z})$ for all $m$, in particular,  determinant of this matrix is $1$. Observe that $w(x,y)$ and $w(x^2,y)$ to not commute in the free group, and hence they do not commute in $SL(2,\Z)$. If $q$ is large enough, their images under the quotient map do not commute in $SL(2,F_q)$. Therefore, either $v_{2,1}$ or $v_{1,2}$ for one of these two matrices is not equal to $0$, and we know $v_{1,1}v_{2,2}-v_{1,2} v_{2,1}=1$ in $F_q$.
We conclude, that for some point in the image of $f$ over $F_q$ either $v_{2,1} (v_{1,1}v_{2,2}-v_{1,2} v_{2,1}) \ne 0$ or  $v_{1,2} (v_{1,1}v_{2,2}-v_{1,2} v_{2,1})\ne 0$. Without loss of generality we can suppose that there exists
a point in the image of $f$ such that $v_{2,1} (v_{1,1}v_{2,2}-v_{1,2} v_{2,1})$.
 
 In this case, we know that there exist at least one solution of the system of the  polynomial equations in $\bar{F_q}$, such that $x_{2,1} (x_{1,1}x_{2,2}-x_{1,2} x_{2,1}) \ne 0$.  Consider a field generated by elements of this solution $x_{1,1}$, $x_{1,2}$, $x_{2,1}$, $x_{2,2}$. This field is clearly a finite
extension of $K$, which we denote by $\mathcal{K}$.
 
 Observe that if $x_{i_1,i_2}$, $i_1, i_2=1,2$ is the solution of the 
system of the equations above over $\mathcal{K}$, then there exist $m$ such that $x_i$ is $m$ periodic point in the group of two times two invertible matrices over $\mathcal{K}$, for polynomial mapping
corresponding to $H$.
Indeed observe that
$$
H_{j_1, j_2}^{(4)^{(2)}}(x_{i_1, i_2},y_{i_1, i_2})=H_{j_1, j_2}^{(4)}(H_{j_1, j_2}^{(4)}(x_{i_1, i_2},y_{i_1, i_2}), y_i) = H_{j_1, j_2}^{(4)}(x_{i_1, i_2}^Q, y_{i_1, i_2}) = (x_{j_1, j_2}^{Q^2})
$$
and, arguing by induction, we obtain that 
$$
H_{j_1, j_2}^{(4^l)}(x_{i_1, i_2},y_{i_1, i_2})=x_{j_1, j_2}^{(q^l)}.
$$

Observe that there exist $l$ such that $x_{i_1, i_2}^{(q^l)}=x_i$ (for $i_1, i_2 = 1, 2$).

Now consider 
$$
\mathcal{K'}= \mathcal{K} [\sqrt{\det x}] =  \mathcal{K} [\sqrt{ (x_{1,1} x_{2,2} - x_{1,2} x_{2,1}} )].
$$
  Recall that we know that $(x_{1,1} x_{2,2} - x_{1,2} x_{2,1}) \ne 0$.
Put $x'=x/ \sqrt{ (x_{1,1} x_{2,2} - x_{1,2} x_{2,1})}$, $x' \in SL(2, \mathcal{K}')$.

Note that  $x',y$ is $4^l$ periodic in $SL(2, \mathcal{K'})$: $x' \ne e$ in  $SL(2, \mathcal{K'})$ and 
the $n$-th iteration $w_{\circ n}$ of $w$ satisfies
$$
w_{\circ n}(x',y)=w(w(\dots w(x',y), y , \dots, y) =x'.
$$

From Remark \ref{obviousremark} we know that it is sufficient to consider words $w$ such that the total number of $x$ in $w$ is equal to zero (otherwise, $w$ is not an iterated identity in some finite cyclic group).
So we assume that $w(x,y)$ is such that the total number of $x$ is equal to zero. Then   $w(e,z) =e$ for all $z$. Therefore, for any periodic point $x'\ne e$ it holds
$$
w_{\circ n}(x',y) \ne e
$$
for some positive integer $n$, and hence $w$ is not an iterated identity for $SL(2, \mathcal{K'})$.

\begin{rem} Alternatively, the first claim of the theorem can be proved by combining the result of Borisov and Sapir about residual finiteness of  the mapping tori (Theorem $1.2$ in \cite{borisovsapir}, rather
than its proof , as explained above) with  characterization of residual finiteness of HNN extension in terms of {\it "compatible"} subgroups
(Theorem $1$ in \cite{moldavanskii}), in case when the corresponding endomorphism is injective, and then reduce the general case in our theorem (when the endomorphism is not necessary injective) to this one.
\end{rem}

\section{Explicit estimates for  Theorem $3.2$  of Borisov and Sapir in \cite{borisovsapir} and the proof of the second part of the Theorem.}\label{se:secondpart}

Theorem \ref{thm:borisovsapir} below is a version of Theorem $3.2$ in \cite{borisovsapir}, which, given an upper bound on the degree of polynomial $D$,
 provides an explicit estimate $Q$ in such a way that if $D$ is equal to $0$ on all solutions of the  system of polynomial equations, than $D$ is equal to zero on the image of the 
 polynomial mapping.
 
For a prime $q$, $F_q$ denotes the filed on $q$ elements.

\begin{thm}\label{thm:borisovsapir}
 Let $q$ be a prime number, $d,n\ge 1$. 
Let $f=f_1, \dots, f_n$ be polynomials on $n$ variables of degree $\le d$, with coefficients in $F_q$, such that $f_i(0, \dots, 0)=0$ 
 for all $i: 1\le i \le n$.
Assume that  $Q$ is a power of $q$ and $D_0\ge 1$ satisfy
$$
Q/D_0>n(n+1)d^{n^2+1}.
$$
Consider a polynomial $D$ of degree at most $D_0$  over $F_q$ on $n$ variables, such that  
$$
D(x_1, \dots, x_n)=0
$$ 
for all $x_i \in \bar{F_q}$ that are solution of the system of the equations 
\begin{equation} \label{fiQ}
f_i(x_1, \dots, x_n) = x_i^Q,
\end{equation}
for all  $i: 1 \le i \le n$.
 Then $D$ is equal to zero on all points in the image of $\bar{F}_q^n$ under $f^{(n)}=\left(f_1^{(n)}, \dots, f_n^{(n)} \right)$.
\end{thm}

We need this theorem in a particular case  when there is no non-zero solution for the system of the equations. In this case it is sufficient to consider $D$ of degree $1$, $D(x_1, \dots, x_n)=x_i$ for some $i$, and we get

\begin{corlem}
Let $q$ be a prime number, $d,n \ge 1$.
Let $f=f_1, \dots, f_n$ be polynomials on $n$ variables of degree $\le D_f$, with coefficients in $F_q$, such that $f_i(0, \dots, 0)=0$ 
 for all $i: 1\le i \le n$. Suppose that there exists $v_1, \dots, v_n \in F_q$ and $i: 1 \le i \le n$ such that $f_i^{(n)}(v_1, \dots, v_n) \ne 0$.
Suppose that $Q$ is a power of $q$ such that   
$$
Q>n(n+1)d^{n^2+1}.
$$

Then the system of equations
$$
f_i(x_1, \dots, x_n) = x_i^Q,
$$
has at least one non-zero solution in $\bar{F_q}$.
\end{corlem}

More precisely, for the proof of Theorem \ref{thm:noid} we need the to find a non-zero solution of the system of $n$ equations, $n=4$, satisfying additionally an inequality $x_1x_4-x_2x_3 \ne 0$, and to obtain such solution
we apply Theorem \ref{thm:borisovsapir} for the polynomials $D(x_1, x_2, x_3, x_4)$ of degree $3$  the form 
$$
D(x_1, x_2, x_3, x_4)=(x_2 x_1x_4-x_2x_3)x_2
$$ 
and 
$$
D(x_1, x_2, x_3, x_4)=(x_2 x_1x_4-x_2x_3)x_3.
$$
(In our matrix notation of the previous section these $x_i$ correspond to  $x_1=x_{1,1}$, $x_2=x_{1,2}$, $x_3=x_{2,1}$, $x_4=x_2,2$).
 For a more general version in Theorem \ref{thm:borsapmol}, we will need to find a system of $4s$ equations, each solution in not proportional to an identity matrix, and the determinants of the corresponding matrices are not equal to zero. To to this, we will apply theorem
\ref{thm:borisovsapir} to the polynomial of degree $3s$, which as a product of the polynomials as above.

Given $(\alpha_1, \dots, \alpha_n)$, $(\beta_1, \dots, \beta_n) \in \Z^n$ we say that $(\alpha_1, \dots, \alpha_n)$ is {\it greater} than  $(\beta_1, \dots, \beta_n)$ in {\it the  lexicographic order} if for the minimal $i$ such that
$\alpha_i-\beta_i \ne 0)$ it holds $\alpha_i>\beta_i$.

Now we recall the definition of the {\it graded lexicographic order} (or for short {\it graded lex order}). Given two monomials $x_1^{\alpha_1} \cdots x_n^{\alpha_n}$ and $x_1^{\beta_1} \cdots x_n^{\beta_n}$, we say that
$x_1^{\alpha_1} \cdots x_n^{\alpha_n}$ is greater than $x_1^{\beta_1} \cdots x_n^{\beta_n}$ in the graded lex order, if either the degree of the first monomial is greater, that is, $\alpha_1 + \cdots + \alpha_n >
\beta_1+\cdots + \beta_n$, or if the degrees are equal (   $\alpha_1 + \cdots + \alpha_n >
\beta_1+\cdots + \beta_n$  ) and $\alpha_1,\dots, \alpha_n$ is greater than $\beta_1, \dots, \beta_n$ in the lexicographic order.

Lexicographic order if a particular case of {\it monomial order}, that is, it is a total ordering of $\mathbb{Z}^n$, satisfying $\alpha+\gamma$ is greater then $\beta+\gamma$ whenever $\alpha$ is greater then $\beta$ and it is a {well-ordering}, meaning that any non-empty subset of $\mathbb{Z}^d$ has a minimal element with respect to this order (see Section 2, Chapter 2 in \cite{coxlittleoshea}).
Fixing a monomial order and given a polynomial $\phi= \sum_i a_{\alpha_{1,i}, \dots, \alpha_{j_i}} x_1^{\alpha_{1,i}} x_2^{\alpha_{2,i}} \cdots x_{n}^{\alpha_{n,i}}$, one can speak about its leading monomial 
$x_1^{\alpha_{1,i}} x_2^{\alpha_{2,i}} \cdots x_{n}^{\alpha_{n,i}}$
 denoted by $LM(\phi)$, and its {\it leading term} $a_{\alpha_{1,i}, \dots, \alpha_{j_i}} x_1^{\alpha_{1,i}} x_2^{\alpha_{2,i}} \cdots x_{n}^{\alpha_{n,i}}$, denoted by $LT(\phi)$.

 This allows to use {\it division algorithm} in $\mathcal{K}[x_1, \dots, x_n]$, $\mathcal{K}$ is some field (see Theorem 3 and its proof, Section 3, Chapter 2 in \cite{coxlittleoshea}). Given an ordered tuple 
$f_1$, \dots, $f_s$ with respect to a fixed monomial order, and given a polynomial $f \in \mathcal{K}[x_1, \dots, x_n]$, the division algorithm procedes as follows. Given $f$, it looks for a minimal $i$ such that
the leading term of $f$ is divided by the leading term of $f_i$, and replaces $f$ by $f- f_i g$, where $g$ is the monomial such that $LT(f) = LT(f_i)g$. At the end we obtain
$$
f=a_1 f_1 + \cdots + a_s f_s +r,
$$
where $r$ is such that no term in $r$ is divisible by a leading term of some $f_i$.

In general, given some tuple $f_s$, this remaining term $r$ is not defined uniquely by the decomposition above, it is difficult therefore to work with the division algorithms. 
This problem no longer occurs if we assume that $f_i$ form a   {\it a Gr{\"o}bner basis} (also called {\it a standard basis}).
There are several equivalent ways to define  a Gr{\"o}bner basis, one 
of such ways is as follows. Elements of $\phi_1$, $f_2$, \dots, $\phi_m \in  \mathcal{K}[x_1, \dots, x_n]$  are said to be  a Gr{\"o}bner basis, if the ideal of  $\mathcal{K}[x_1, \dots, x_n]$ generated by
leading terms of $\phi_j$ is equal to the ideal, generated by the leading term of the ideal $I$ generated by $\phi_j$ (Definition 5 in \cite{coxlittleoshea}).
 In general, the ideal of leading terms of $I$ can be larger than that generated by leading terms of $\phi_j$. For a Gr{\"o}bner basis this can not happen, and this provides an effective way to determine whether a polynomial
$F$ belongs to the ideal  $I$ generated by $f_j$.

For example, given polynomials $f_i$ ($i: 1 \le i \le n$) on $x_1$, $x_2$, \dots, $x_n$, consider polynomials $\phi_i=f_i -x_i^Q$. Suppose that that the degrees of $f_i$ are smaller than $Q$, for all $i: 1 \le i \le n$.
Then the leading terms of $f_i$ with respect to graded lex order is $x_i^Q$. The Least Commond Multiple of the Leading Monomials of $\phi_i$ and $\phi_j$ is equal to their product , $x_i^Qx_j^Q$, and hence
by Diamond Lemma $\phi_i$ is a Gr{\"o}bner  basis (see e.g. Theorem 3 and Proposition 4 in Ch.2, Sect.9 of \cite{coxlittleoshea}).

A straightforward observation in  Lemma  \ref{le: lemma3.3} below is Lemma $3.3$ from \cite{borisovsapir}. While the formulation of that lemma in \cite{borisovsapir} states "for $Q$ large enough", the proof shows that it is sufficient to take $Q$ which greater than the maximum of the  degrees of $f_i$, as stated below, and it is not difficult to estimate this codimension.

\begin{lem} \label{le: lemma3.3} [a version of  Lemma 3.3 in \cite{borisovsapir}]  \label{le:codimension}

Let $f_1$, \dots, $f_n$ be polynomials on $x_1$, \dots, $x_n$ over $F_q$, $q$ is a prime number. Take $Q$ such that $Q$ is greater than the maximum of the degrees of $f_i$, $i: 1 \le i \le n$.
Let $I_Q$ be the ideal in $\bar{F}_q[x_1, \dots, x_n]$ generated by polynomials $f_i(x_1, \dots, x_n)- x_i^Q$, $i: 1 \le i \le n$. Then $I_Q$ has finite codimension in $\bar{F}_q[x_1, \dots, x_n]$, this codimension is at most
$Q^n$.

\end{lem}

This shows in particular, that any solution (in $\bar{F}_q$ the system of equations $f_i(x_1, \dots, x_n) = x_i^Q$  belongs to a finite extension of $F_q$, of degree at most
$Q^n$.
 Observe that a cardinality of such finite extension is at most $q^{Q^n}$.

{\bf Proof.}  
Observe that if a monomial on $x_1$, $x_2$, \dots , $x_n$ is divisible by $x_i^Q$, then this monomial is equivalent $\pmod {I_Q}$ to a sum of monomials of lesser degree.
Indeed, 
$$
x_i^Q \prod x_i^{\alpha_i}  \equiv f_i(x_1, \dots, x_n )   \prod x_i^{\alpha_i}.
$$
 All monomials of $f_i$ have degree strictly lesser than $Q$, and the degree of monomials on the right hand side is therefore
lesser than $Q + \sum_i \alpha_i$.
Therefore, any polynomial on $x_1$, $x_2$, \dots , $x_n$ is equivalent  $\pmod {I_Q}$ to a linear combinations of monomials of the form $\prod x_i^{\alpha_i}$, such that $\alpha_i < Q$ for all $i$. Observe that
the number of such monomials is $Q^n$

\begin{rem} We will use only the (trivial) upper bound for the codimension, but it is not difficult to see that
the codimension is in fact equal to $Q^n$. Indeed, as we have already mentioned the   $f_i-x_i^Q$ is a Gr{\"o}bner basis, and by a "Diamond lemma"  we know that any   linear combination of $\prod x_i^{\alpha_i}$,  with at least one
non-zero coefficient ,  such that $\alpha_i < Q$ does not belong to $I_Q$.

\end{rem}

Given polynomial $f_i$, $1 \le i \le n$ over some field $\mathcal{K}$, we can consider a mapping $f = (f_1, f_2, \dots, f_n)$ from $\mathcal{K}^n$ to $\mathcal{K}^n$. For $j\ge 1$ we denote by
$f^{(j)}= (f_1^{(j)}, \dots, f_n^{(j)})$ its $j$-th iteration.

We recall in Lemma \ref{le:lem3.4} below Lemma 3.4 in \cite{borisovsapir}), for the convenience of the reader we recall its proof.

\begin{lem} \label{le:lem3.4} [Lemma $3.4$ in \cite{borisovsapir}]
Let $f_i$ be polynomials on $x_1$, \dots, $x_n$ with coefficients in $F_q$, $q$ is a prime number, and $Q$ be a power of $q$.
For each $j\ge 1$ the $j$-th iteration of the polynomial mapping $f=\left( f_1, \dots, f_n\right)$ satisfies for all $i: 1 \le i \le n$
$$
f_i^{(j)} - x_i^{Q^j} \in I_Q,
$$
where $I_Q$ is the ideal generated by $f_i-X_i^Q$, $i: 1 \le i \le n$.

\end{lem}

{\bf Proof.}
The proof is by induction on $j$. Suppose that the statement is true for all $j\le m$.
Observe that 
$$
f_i^{(m+1)}(x_1, x_2, \dots, x_n) = f_i \left(f_1^{(m)}, f_2^{(m)}, \dots, f_n^{(m)}     \right) \equiv f_i \left(x_1^{Q^m}, \dots, x_n^{Q^m}   \right) 
$$
The last congurence above  follow from the induction hypothesis for $j=m$.
Observe also that  since $Q$ is a power of $q$, over any field of characteristic $q$ it holds
$$
f_i \left(x_1^{Q^j}, \dots, x_n^{Q_j}   \right) = f_i(x_1, \dots, x_n)^{Q^j} \equiv x_i^{Q^{j+1}} \pmod  I_Q,
$$
the last congruence is a consequence of the induction hypothesis for $j=1$.

\begin{lem} \label{le:algebraicdependance} Let $F_1$, $F_2$, $F_{n+1}$ are polynomials on $x_1$, $x_2$, \dots, $x_{n}$ over some field $\mathcal{K}$. Suppose that the degrees of $F_i$ are $\le d$.
If $s>0$  is such that the binomial coefficients satisfy
$$
C_{s+n+1}^{n+1} \ge C_{sd+n}^n,
$$
then there exists a  non-zero polynomial $\Psi$  over $\mathcal{K}$ on $n+1$ variables of degree at most $s$ such that
$$
\Psi(F_1, \dots, F_{n+1})=0.
$$
The assumption on $s$ is in particular satisfied if 
$$
s \ge (n+1)d^n
$$

\end{lem}

In the lemma above,  it is essential that $s \ge {\rm Const} \cdot d^n$.

\medskip

{\bf Proof.}  It is clear that it is sufficient to consider the case when at least one of $f_i$ has at least one non-zero coefficient.
Take some integer $s$ and a polynomial $\Psi$ of degree $s$. Let us  compute the number of possible monomials  on $n$ variables $x_1$, \dots, $x_n$ in $\Psi(F_1, \dots, F_{n+1})$. All monomials are of degree at most
$s d$, that is , of the form   $X_1^{\beta_1} X_2^{\beta_2} \cdots X_{n}^{\beta_{n}}$, $\beta_i \ge 0$, $\sum_j\beta_j \le sd$. This is the number to write $sd$ as the sum of $n+1$ non-negative summands, 
which is equal to
 $C_{sd+n}^n$.
Consider possible monomials of degree $\le s$  on $n+1$ variables, they are of the form 
 $y_1^{\alpha_1} y_2^{\alpha_2} \cdots y_{n+1}^{\alpha_{n+1}}$, $\alpha_i \ge 0$, $\sum_j\alpha_j \le n$  and hence their number  is equal to $C_{s+n+1}^{n+1}$.

Take $s$ such that 
$$
C_{s+n+1}^{n+1} \ge C_{sd+n}^n.
$$
Observe that  there exists a non-zero polynomial $\Psi$ of degree at  most $s$ such that 
$$
\Psi(F_1, \dots, F_{n+1})=0.
$$
 Indeed, if we consider the coefficients of the polynomial $\Psi$  (taking value in the field $\mathcal{K}$) as variables, we get  $C_{sd}^n$ linear equations  on at least  $C_{s+n+1}^{n+1}$ variables.
Since the number of variables greater or equal to the number of linear equation,  this  system has at least one non-zero solution over $k$.

Finally, observe that
the assumption on $s$ in the formulation of the Lemma   is satisfied if 
$$
(s+1)  (s+2) \cdots  (s+n+1) \ge (n+1) (sd+1) \cdots (sd+n), 
$$
and the latter is satisfies whenever $s+n+1 \ge s \ge (n+1)d^n$.

Lemma \ref{le:algebraicdependance} allows us to obtain explicit estimates for Lemma $3.5$ in \cite{borisovsapir}:

\begin{lem}  \label{le:l3.5}
Given $d$,
take an integer $Q$ which is a power of a prime  $q$ such 
that $Q > (n+1) d^{n^2}$ and $k=(n+1) d^{n^2} $.
Consider polynomials $f_1$, $f_2$, \dots, $f_n$ over $F_q$ on $n$ variables.  Suppose that the degrees of $f_i$ are $\le d$.
Let $a_1$, \dots, $a_n$  in the algebraic closure  $\bar{F}_q$ of $F_q$ are the solution of the system of equations
$$
f_i(a_1, a_2, \dots, a_n) = a_i^Q.
$$
Then for all $i: 1 \le i \le n$ the polynomial
$$
(f_i^{(n)}(x_1, \dots, x_n) -f_i^{(n)}(a_1, \dots, a_n))^k
$$
is contained in the localization of $I_Q$ at $a_1, \dots, a_n$.
As before, $I_Q$ denotes the ideal in $\bar{F}_q[x_1, \dots, x_n]$ generated by polynomials $f_i(x_1, \dots, x_n)- x_i^Q$, $i: 1 \le i \le n$.
\end{lem}

{\bf Proof.}
For each $i: 1 \le i \le n$ consider the $i$-th coordinate of the  iterations of $f$: $F_{1,i} =x_i$, $F_{2,i}(x_1, \dots, x_n)=f_i(x_1, \dots, x_n)$, $F_{3,i}(x_1, \dots, x_n) = f_i^{(2)(x_1, \dots, x_n)}$, 
\dots, $F_{n+1,i}(x_1, \dots, x_n)=f_i^{(n)(x_1, \dots, x_n)}$.
Observe that for all $j$ the degree of $F_{j,i}$ is at most $d^n$.

Apply lemma \ref{le:algebraicdependance} to $F_{1,i}$, $F_{2,i}$, \dots, $F_{n+1,i}$. We conclude that for each $i: 1 \le i \le n$ there exists a non-zero polynomial $\Psi_i$ over $F_q$ on $n+1$ variables of degree at most $(n+1)(d^n)^n$ such that
$$
\Psi_i(x_i, f_i(x_1, \dots, x_n), f_i(x_1, \dots, x_n), \dots, f_i(x_1, \dots, x_n))=0
$$

The rest of the proof follows the argument from \cite{borisovsapir}:
using the fact that $f_i^{(j)} - x_i^{Q^j} \in I_Q$ (see Lemma  \ref{le:lem3.4}), we can rewrite
$$
\Psi_i(x_i, f_i, f_i^{(2)}, \dots, f_i^{(n)})
$$
 as a polynomial $P_{Q,i}$ in one variable  $x_i$ modulo $I_Q$. Since
  $\Psi_i(x_i, f_i, f_i^{(2)}, \dots, f_i^{(n)})=0$,  we have  $P_{Q,i}(x_i) \in I_Q$

By the assumption of the lemma,  $Q> (n+1)d^{n^2}$,  and hence $Q$ is larger than the degree of $\Psi_i$. Observe that in this case the polynomial in $x_i$  we get is not
 zero. (Indeed, take maximal $j$ such that  $y_j$ is present at least in one monomial of $\Psi$;
 among  monomials of $\Psi$ consider those  where the degree of $y_j$ is maximal. Among such monomials, if there several like this, take maximal $j'$ such that $y_{j'}$ is present, take a monomial where its degree is maximal, etc. In this way we obtain some monomial
in $\Psi$ which will give maximal degree of $x_i$ for $P_{Q,i}$).
Note that the degree of $P_{Q,i}$ is at most $Q^n \deg \Psi \le Q^n (n+1) d^{n^2}$.

Write
$$
P_{Q,i}(x_i) = \sum_{m=1}^M b_m  (x_i-a_i)^m,
$$
here $b_m \in \bar{F}_q$, $1 \le m \le M$ are such that  $b_M \ne 0$.

It is clear that  $M \le  \deg P_{Q,i} \le Q^n (n+1) d^{n^2}$, and in particular 
$$
P_{Q,i}(x) = (x-a_i)^L u(x),
$$
 where $L \le M \le Q^n (n+1) d^{n^2}$  and the polynomial $u(x)$ is such that $u(a_i) \ne 0$.
Recall that by the assumption of the lemma $k =  (n+1) d^{n^2}$. It is essential for the proof that $k$ does not depend on $Q$.

Since $P_{Q,i}(x_i) \in I_Q$, 
we conclude that $(x_i-a_i)^L \in I_Q^{a_1, a_2, \dots, a_n}$.
We have $Q^nk \ge L$, and therefore $(x_i-a_i)^{Q^n k} \in I_Q^{a_1, a_2, \dots, a_n}$.

By the assumption of the Lemma, $f_i(a_1, \dots, a_n)=a_i^Q$. Since the characteristic of the field is $p$, this  implies that $f_i^{(m)}(a_1, \dots, a_n) = a_i^{Q^m}$ for all $m\ge 1$ . Hence by Lemma \ref{le:lem3.4} we obtain
$$
 f_i^{(n)}(x_1, \dots, x_n) - f_i^{(n)}(a_1, \dots, a_n)= f_i^{(n)}(x_1, \dots, x_n) - a_i^{Q^n} \equiv x_i^{Q^n} - a_i^{Q^n} \pmod {I_Q^{(a_1, \dots, a_n}}
$$
Since the characteristic of the field is $p$ and $Q$ is a power of $p$, we know that  $x_i^{Q^n} - a_i^{Q^n}=(x_i-a_i)^{Q^n}$.
Therefore we can conclude that 
$$
\left( f_i^{(n)}(x_1, \dots, x_n) - f_i^{(n)}(a_1, \dots, a_n) \right)^k  = (x_i-a_i)^{kQ^n} \equiv 0 \pmod {I_Q^{a_1, \dots, a_n}}
$$

As a corollary, we obtain a version of Lemma $3.6$ in \cite{borisovsapir}.

\begin{corlem} \label{le:3.6}
Let $q$ be a prime number, $d,n\ge 1$.
Consider $n$ polynomials $f_i$, $1 \le i\le n$  on $n$ variables, with coefficients in $F_q$,   of degree at most $d$.
Take a polynomial $D$ with  coefficients  in $F_q$,
which vanishes on all solutions 
 in $\bar{F_q}$ of the system of the equations
$$
f_i(a_1, a_2, \dots, a_n) = a_i^Q
$$
Assume, as in  Lemma \ref{le:l3.5}, that  $Q > (n+1) d^{n^2}$ and  $k= (n+1) d^{n^2} $. Put $K= (k-1)n+1$.
Then for any $a_i$, $1\le i \le n$ which the solution of the above mentioned system of polynomial equations
$$
\left(D(f_1^{(n)}(x_1, x_2, \dots, x_n),..., f_n^{(n)}(x_1, x_2, \dots, x_n) ) \right)^K = 0 \pmod {I_Q^{(a_1, \dots, a_n)}}.
$$
We recall that $I_Q$ denotes the ideal in $\bar{F}_q[x_1, \dots, x_n]$ generated by polynomials $f_i(x_1, \dots, x_n)- x_i^Q$, $i: 1 \le i \le n$.

\end{corlem}

{Proof of Corollary  \ref{le:3.6}.}
Take $a_i \in \bar{F}_q$ such that $f_i(a_1, a_2, \dots, a_n) = a_i^Q$. 
We have $f_i^{(j)}(a_1,\dots, a_n) = a_i^{Q^j}$, for all $j\ge 1$.
Rewrite $D(x_1, \dots, x_n)$ as a  polynomial in $x_i-a_i^{Q^n}$, ($1 \le i \le n$), that is, 
$$
D(x_1, x_2, \dots, x_n) =E(x_1 -a_i^{Q^n}, \dots,   x_n- a_n^{Q^n}),
$$
 where $E$ is a polynomial (depending on $a_1$, \dots, $a_n$) with coefficients in $\bar{F_q}$.
 Since 
$$
f_i(a_1, a_2, \dots, a_n) = a_i^Q
$$ for all $i$, we know  by the assumption of the Corollary that
$D(a_1, \dots, a_n)=0$. Hence 
$$
E(0,0,\dots, 0) =D(a_1^{Q^n}, \dots, a_n^{Q^n})= D(a_1, \dots, a_n)^{Q^n}  =0,
$$
and therefore the polynomial $E$ does not have a free term. 
Observe that $D^K$ can be therefore written as sum of monomials in $x_i-a_i^{Q^n}$. Since  $K\ge (k-1)n+1$,
for each of these monomials  there exists $i$, $1\le i \le n$ such that this monomial is divisible by at $(x_i-a_i^{Q^n})^k$.
Therefore, $\left(D(f_1^{(n)}(x_1, x_2, \dots, x_n),..., f_n^{(n)}(x_1, x_2, \dots, x_n)) \right)^K$ is congruent $\pmod {I_Q}$ to a sum of polynomials, 
for each of these polynomial there exists $i: 1\le i \le n$ such that the polynomial is   divisible by $(f_i^{(n)}- a_i^{Q^n})^k$. 

In other words, each of the above mentioned polynomials is divisible by 
$$
(f_i^{(n)}(x_1, \dots, x_n)- f_i^{(n)}(a_1, \dots, a_n))^k.
$$
Applying Lemma \ref{le:l3.5} we conclude  that each of these polynomials belong to ${I_Q^{(a_1, \dots, a_n)}}$, and hence their sum belongs to ${I_Q^{(a_1, \dots, a_n)}}$

{\bf Proof of Theoreom \ref{thm:borisovsapir}.}

By the assumption of the theorem,
$$
Q/D_0>n(n+1)d^{n^2+1}, 
$$
and hence
$$
Q/D_0 > dn((n+1)d^{n_2}-1)+1.
$$
We will prove the theorem under the assumption above.

Observe that $Q >D_0 dn((n+1)d^{n_2}-1)+1 \ge d+1 >d$. This shows that $Q$ is
$Q$ is greater than the degrees of $f_i$. We have already mentioned that in this case we know that $f_i-x_i^Q$  form a Gr{\"o}bner basis with respect to Graded Lex order.
Recall that in this situation no non-zero  polynomial of degree strictly smaller than $Q$  belongs to the ideal generated by  $f_i-x_i^Q$, $1\le i \le n$. 
In particular,  if we assume that the degree of the polynomial 
$$
P=\left(D(f_1^{(n)}(x_1, x_2, \dots, x_n),..., f_n^{(n)}(x_1, x_2, \dots, x_n) )\right)^K
$$
 (where $D$ and $K$ are as in Lemma \ref{le:3.6}), is strictly less then $Q$, we conclude that
$P$ does not belong to the ideal $I_Q$ generated by $f_i-x_i^Q$.

Take a polynomial $D$ satisfying the assumption of Theorem \ref{thm:borisovsapir} which is zero on all the solutions of the system of equations, and non-zero at at least one point of the image of $f$.
We want to obtain a contradiction.

Observe that since 
$$
Q/D_0 > dn((n+1)d^{n_2}-1)+1.
$$
we know  that $Q/D_0>   d ((k-1)n)+1$ for  $k= (n+1) d^{n^2} $, 
and that  $Q > (n+1) d^{n^2}$. 

Put $K= (k-1)n+1 =  ((n+1) d^{n^2} -1) n+1$.

Observe  that  $K$  and $k$ satisfy the assumption of  the corollary  \ref{le:3.6}. It is essential for our argument that  $K$ does not depend on $Q$. Observe that the polynomial $D(f_1, \dots, f_n)$ has at least one non-zero coefficient, since from the assumption of the theorem we know that this polynomial takes at least one non-zero value. This implies that the polynomial
$$
\left(D(f_1^{(n)}(x_1, x_2, \dots, x_n),..., f_n^{(n)}(x_1, x_2, \dots, x_n)) \right)^K.
$$
has at least one non-zero coefficient.

This polynomial above belongs to $I_Q^{(a_1, \dots, a_n)}$ for any solution of the system of equations $a_1, \dots, a_n$, with $Q$ satisfying the assumption of the corollary. Now, like in the second version of
the proof of \cite{borisovsapir}, observe that if $a_1, \dots, a_n$ is not a solution of the system of the equations \ref{fiQ}, then the localisation of $I_Q$ at $a_1, \dots, a_n$ is the whole ring of polynomials 
$\bar{F_q}[x_1, \dots, x_n]$.

Indeed,  we know in this case that there exists $i$, $1 \le i \le n$ such that  
$$
(f_i-x_i^Q)(a_1, \dots, a_n) \ne 0.
$$
 Then  $1/(f_i-x_i^Q)$ belongs to the localisation, and since $f_i-x_i^Q$ belongs to $I_Q$, we conclude that $1 \in I_Q^{a_1, \dots, a_n}$.

We know therefore  that for any $a_1, ..., a_n \in \bar{F_q}^n$  (whether it is a solution of the system of equations of whether it is not)
the polynomial $D(f_1, \dots, f_k)^K$ belongs to the localisation of $I_Q$ at $a_1, \dots, a_n$, and hence
$D(f_1, \dots, f_k)^K$ belongs to $I_Q$. 
But in this case the degree of $D(f_1, \dots, f_k)^K$ is greater of equal to $Q$.

This completes the proof of Theorem \ref{thm:borisovsapir}.

\medskip

{\bf Proof of Theorem \ref{thm:noid}.}
First we observe again that it is sufficient to consider words on two letters:

\begin{rem} Let $m\ge 2$.
For each $m$ fix  $u_2(x,y)$, \dots, $u_m(x,y)$ in the free group generated by $x$ and $y$  such $x_1$, $u_2(x,y)$, \dots, $u_m(x,y)$   freely generate a free subgroup on $m$ generators.
 For any word 
$\bar{w}(x_1, \dots, x_m)$, not freely equivalent to an empty word,
consider the following word on two letters:  
$$
w(x_1, x_2) = \bar{w}(x_1, u_2(x_1, x_2),u_m(x_1,  x_2)).
$$
 As we have mentioned already, this word is not freely equivalent to an empty word, and $\bar{w}$ is an iterated identity in some group whenever this is the case for $w$.
Now assume in addition that for each $j: 2 \le j \le m$ there is a single occurance of $x$ or $x^{-1}$ in the word $u_j(x,y)$. (For example, one can take $u_j(x,y)=  y^j x y^{-j}$).
Then the  total number of $x_1$ and $x_1^{-1}$ in $w$ and the length of  $\bar{w}$ satisfy $l_{x_1}(w) \le  l(\bar{w})$

\end{rem}

\begin{rem} 

Take a word $w(x,y)$ of length $l$. Then $w_{\circ 4}$ has length at most $l^4$. 
The four polynomials  in $x_{1,1}$, $x_{1,2}$, $x_{2,1}$, $x_{2,2}$ for the entries of $w(x,y)$ have degree at most $l$
The polynomials for  for the entries of  $w_{\circ 4}(x,y)$ have degree st most  $l^4$.

\end{rem}

\begin{rem} \label{rem:twomatrices}
Take
\[ \bar{x}= \left( \begin{array}{ccc}
1 & 2  \\
0 & 1
 \end{array} \right) 
\mbox{ and  }
 y= \left( \begin{array}{ccc}
1 & 0  \\
2 & 1
 \end{array} \right)\] 
 Take a product of $L$ terms of the form $\bar{x}$, $\bar{x}^{-1}$, $y$ or $y^{-1}$. Then 
the entries of the resulting matrix are at most $3^L=\exp(\ln3 L)$.
For a product of $L$ terms of the form $\bar{x}^2$, $\bar{x}^{-2}$, $y$ or $y^{-1}$ the entries of the resulting matrix are at most $6^L=\exp(2\ln 3 L)$.
\end{rem}

Let $w$ be a word on $x_1$ and $x_2$ of length at most $l$. Take $y$ as in Remark \ref{rem:twomatrices} (this $y$ satisfies Convention \ref{strongerconv:freesubgroup}),
consider rationals functions $R_i$ 
in $x_1$, \dots, $x_4$ which are entries for $w(x,y)$, and the corresponding polynomials $H_i$. For each $j$ it holds 
$$
R_{j_1,j_2}(x_{1,1}, x_{1,2}, x_{2,1}, x_{2,2}) = H_{j}(x_{1,1}, x_{1,2}, x_{2,1}, x_{2,2}) /(x_{1,1}x_{2,2}-x_{1,2}x_{2,1})^s,
$$
 where $s$ is the number of occurrences  of $x^{-1}$ in $w$.
Observe that the coefficients of these polynomials $H_{i,j}$  satisfy the assumption of Remark \ref{rem:twomatrices}, and hence their coefficients 
are at most $4^l$.  

Take $x$ as in Remark \ref{rem:twomatrices}.
Since $\bar{x}$ and $y$ generate a free group on two generators and since $w(x,y)$ is not freely equivalent to an empty word, we know 
 that $w(\bar{x},y)$ is not an identity matrix. Moreover, we know that $w(x,y)$ does not commute with $w(x^2,y)$ in the free group generated by $x$ and  $y$, and hence
 the matricies $w(\bar{x},y)$  $w(\bar{x}^2,y)$ do not commute, that is their commutator $[w(\bar{x},y),w(\bar{x}^2,y)]$ is not an identity matrix.

  The coefficients of $w_{\circ 4}(\bar{x},y)$ are at most $\exp(Cl^{4})$ and the coefficients of
$w_{\circ 4}(\bar{x}^2,y)$ are at most $\exp(2Cl^4)$, for $C=\ln 3$. Since these matrices do not commute,  for at least one of these matrices either $x_2 \ne e$ or $x_3 \ne e$.
We conclude that there exist intergers $x_i \in \mathbb {Z}$, in the image of $w^{(4)}$,  such that the matrix $x$ they form is in $SL(2,\Z)$, such that their coefficients are at most $\exp(2Cl^4)$, for $C=ln 3$ and such that the matrix
$x$ is not a diagonal matrix (that is, either $x_2 \ne e$ or $x_3 \ne e$). We want to find a prime number $q$ , such that the image of the matrix $x$ over quotient map to $F_q$ is not a diagonal matrix. That is, we want that 
 the coefficients of the above mentioned
 matrix  $x$ modulo $q$ satisfy $x_2 \ne e$ or $x_3 \ne e$ in $F_q$. 

Recall that Prime Number theorem says that the number $\phi(x)$ of prime numbers smaller than $x$ satisfies $\phi(x) /(x/\ln(x)) \to 1$ as $x \to \infty$. In particular, we know that
the number of primes between $x/2$ and $x$ is $(1\pm \epsilon)x/(2 \ln x)$, for any $\epsilon>0$ and all sufficiently large $x$.
This implies that  for  any positive integer  $M$,
there is a prime $q$, $q \le C' \ln{M}$,  
 such that $M$ is not divided by $q$. If $M$ is sufficiently large, we can take $C'$  to be close to $1$.

We can therefore choose a prime  $q \le C' Cl^4$, with $C=2\ln(3)$ and $C'$ close to one if $l$ is large enough,   such there is a non-diagonal  matrix  in $SL(2, F_q)$ in the image of the polynomial mapping corresponding to $w^{(4)}$. In this case, $f_i$ and $v$, considered over $F_q$ satisfy the assumption on the second part of Theorem \ref{thm:borisovsapir}.
Consider  the polynomial  $D(x_{1,1}, x_{1,2}, x_{2,1}, x_{2,2}) = x_{1,2}(x_{1,1}x_{2,2}-x_{1,2}x_{2,1})$, and the polynomial $ D_2(x_{1,1}, x_{1,2}, x_{2,1}, x_{2,2}) =   x_{2,1}(x_{1,1}x_{2,2}-x_{1,2}x_{2,1})$. The degree of these two polynomials is equal to $3$. Observe that there exists a point $v_{1,1}, v_{1,2}, v_{2,1}, v_{2,2}$ in the
image of $F_q^4$ of the polynomial mappings corresponding to  $H^{(4)}$, such that either $D(v_{1,1}, v_{1,2}, v_{2,1}, v_{2,2})\ne 0$ or $D_2(v_{1,1}, v_{1,2}, v_{2,1}, v_{2,2})\ne 0$. Without loss of the generality we can assume that
$D(v_{1,1}, v_{1,2}, v_{2,1}, v_{2,2})\ne 0$

Taking $Q$ satisfying the assumption of  Theorem \ref{thm:borisovsapir} for $D_0=3$, $d=l$,  $n=4$, that is
$$
Q> 3 \times 20 \times l^{17}.
$$
we conclude that the system of equations $H_{i_1, i_2}(x_{1,1}, x_{1,2}, x_{2,1} x_{2,2}) = x_{i_1, i_2}^Q$ has a solution over the algebraic closure of $F_q$, such that  $D(x_{1,1}, x_{1,2}, x_{2,1}, x_{2,2})\ne 0$.
If this is the case, by Lemma $3.1$ we know that the solution belongs to a finite extension  $\mathcal{K}$ of $F_q$, of degree at most $Q^4$.
The number of elements in this field is $q^{Q^4}$.

Consider $\mathcal{K'}= K(\sqrt{x_{1,1}x_{2,2}-x_{1,2}x_{2,1}})$. The cardinality of $\mathcal{K'}$ is at most $q^{2Q^4}$. Taking in mind that we can chose $q \le C' Cl^4$, with $C=2\ln(3)$ and $C'$ close to one if $l$ is large enough
and $Q = 61 l^{17}$, we see that we can chose the field  $K'$ as above of cardinality at most $(2l^{4})^{61^4 l^{17 \times 4}} \le \exp(l^{68+\epsilon})$ for all $l\ge L$.

Since $D(x_{1,1}, x_{1,2}, x_{2,1}, x_{2,2})\ne 0$,  we know that
$x_{1,1}x_{2,2} - x_{1,2} x_{2,1} \ne 0$,  and $x_{1,2}\ne 0$. Dividing $x$ by $\sqrt{x_{1,1}x_{2,2} - x_{1,2} x_{2,1}}$ we obtain a non-identity solution in $SL(2, \mathcal{K'})$ of the system of the equations
$$
R_{i_1, i_2}= x_{i_1, i_2}^Q.
$$

Observe that the cardinality of $SL(2)$ over a finite field of cardinality $N$ is $N^3-N \le N^3$.
Therefore,  the number of elements of
$SL(2, \mathcal{K'})$  is at most $\exp(l^{68+\epsilon})$, for all $l> L$.

As in the previous section, we observe that any non-identity solution $x$ in $SL(2, \mathbb{Z})$ of the above mentioned system of equations provides a periodic point : $w_{\circ m}(x)=x$ for some $m\ge 1$.
And we can remark again that if $w$ represents an element of a commutator group in the free group, then $x\ne e$, $w_{\circ m(x)}=x$ implies that $w_{\circ m'}(x) \ne e$ for all $m' \ge 1$.

\section{General dynamics $a_1, \dots, a_k \to w_1(a_1, \dots, a_k), \dots, w_k(a_1, \dots, a_k)$. Proof of Theorem  \label{thm:borsapmol}} \label{se:generaldynamics}

In the definition of iterated identities we consider the iteration on the first letter. Now more generally consider  $s\ge 1$ and words $w_1(a_1, \dots, a_s)$, \dots, $w_k(a_1, \dots, a_s)$,  the mapping

$w: a_1, \dots, a_s \to w_1(a_1, \dots, a_s), \dots, w_k(a_1, \dots, a_s)$ and its iterations: 
$$
w_{\circ n}(a_1, \dots, a_n) 
=w_1 (w_{1, \circ n-1}(a_1, \dots, a_s), \dots , w_{k, \circ n-1}(a_1, \dots, a_s)), \dots,
$$
$$
\dots w_{s, \circ n-1}(a_1, \dots, a_s)).
$$
For some tuples of words, in contrast when the iteration only on the first letter is allowed, it may happen that some  iteration of $w$  is freely equivalent to the  identity, in this case its image is trivial in the free group, and hence in any other group.
For example, if $k=2$ and $w_1(a_1,a_2) = w_2(a_1, a_2) = [a_1, a_2]$, then it is clear that  $w_{1,\circ 2}(a_1, a_2) = w_{2, \circ 2}(a_1, a_2) = [[a_1,a_2], [a_1,a_2]] \equiv e$ in the free group generated by $a_1$ and $a_2$.

In fact, it is possible that all the  coordinates of the first $n-1$ iterations are not equal to one in the free group, and all coordinates on the $n$-th iteration is equal to one:

\begin{exa}
Consider words $w_i$ on $x_1, \dots, x_n$, $n\ge 2$:
$w_1= [x_1, x_n]$, $w_2 = [x_1,x_n]$, $w_3=[x_2,x_n]$,  $w_i= [x_{i-1}, x_n]$ for $i \ge 3$.
Then for all $m\ge n-1$ and all $i$ ($1 \le i \le n$) the word $w_{m,i}$ is not freely equivalent to an empty word. For all $m \ge n$ and all $i$ ($1 \le i \le n$)
the word $w_{m,i}$ is not freely equivalent to an empty word.
\end{exa}

{\bf Proof.} Note that 
$$
x_1 \to [x_1, x_n] \to [ [x_1, x_n], [x_{n-1}, x_n]] \to ...
$$
$$
x_2 \to [x_1, x_n] \to [ [x_1, x_n], [x_{n-1}, x_n]] \to ...
$$
$$
x_3 \to [x_2, x_n] \to [ [x_1, x_n], [x_{n-1}, x_n]] \to ...
$$
$$
x_4 \to [x_3, x_n] \to [ [x_2, x_n], [x_{n-1}, x_n]] \to ...
$$

We see that for all $m\ge 1$ the images of $m$-th iteration evaluated at  $x_1$, $x_2$, \dots, $x_{m+1}$ are equal.
In particular, for $m=n-1$ the image of $n-1$-th iteration takes the same value at all $x_i$, and hence $w_{n,i}$ is freely equivalent to an empty word for all $i$ ($1 \le i \le n$).
Observe, that if  for some $k$  the elements  $y_1$, $y_2$, \dots, $y_k$ freely generate a free group on $k$ generators, then $[y_1,y_k]$, \dots, $[y_{k-1}, y_k]$ freely generate a group on $k-1$ generators.
Using this fact and arguing by induction on $j$ we observe for all $j<n$ the elements  $w_{j,i}$, $i : n-j +1 \le  i \le n$ freely generate a group on $n-j$ generators.
This implies in particular that for $j<n$ all coordinates of the $j$-th iteration are non-trivial.

A generalization of the first part of Theorem \ref{thm:noid} says that if the all components of the iteration map are not trivial in a free group, then there is a finite group where where all components of the iterations remain non-trivial.

\begin{rem} Suppose that the words $w_1(x_1, \dots, x_n)$, \dots, $w_n(x_1, \dots, x_n)$ are such that $w_1(x_1, \dots, x_n)$, \dots, $w_n(x_1, \dots, x_n)$ generate a free subgroup of rank $n$ in the free group generated by $x_1$, \dots, $x_n$.
Then for all $m \ge 1$ and all $i$ , $1 \le i \le n$ the iteration $w^i_{\circ m} \ne e$ in the free group generated by $x_1$, \dots, $x_n$.

\end{rem}

{\bf Proof.} Observe that  the endomorphism $w: x_i \to w_i(x_1, \dots, x_n)$ is injective, since otherwise the free group $F_n$ would have a quotient over non-trivial normal subgroup with the image isomorphic
to $F_n$. It is well known and not difficult to see that this can not happen, in other words, the free group (as any other residually finite gorup) is Hopfian (see e.g. Thm 6.1.12 in \cite{robinsonbook}).
Therefore, any itetation  $w_{\circ m} $ of the endomorphism $w$ is injective. Hence the image of  $w^i_{\circ m}$ is isomorphic to $F_n$, that is, this image  is a free group of rank $n$.
This implies in particular that for all $m\ge 1$ and all $i$  $w^i_{\circ m} \ne e$ in the free group generated by $x_1$, \dots, $x_n$.

\begin{rem}
Suppose that  $w^i_{\circ m} \ne e$ in the free group, for all $m\le n$ and all $i$.
 Then $w^i_{\circ m} \ne e$ for all $m$ and all $i$.
\end{rem}

{\bf Proof.} Consider images of the free group $F_n$ (generated by $x_1$, \dots, $x_n$) with respect to $w$, $w_{\circ 2}$, \dots, $w_{\circ n}$. If there exists at list some element , not equal to $e$, in the image of $w_{\circ n}$,
then there exists $m<n$ such that the rank of the free group in the image of $w_{\circ m}$ is equal to that in the image of $w_{\circ m+1}$, and this rank is at least $1$. In this case the restriction of $w$ to the image of $w_{\circ m}$ is injective, that is,
if $g,h$ in  the image of $w_{\circ m}$ are such that $w(g)=w(h)$, then $g=h$. Arguing by induction we see that  for all $t\ge 1$ the restriction of $w_{\circ t}$ to the image of $w_{\circ m}$ is injective.
Therefore, if $w_{\circ m} (x_i) \ne e$ in the free group generated by $x_1$, \dots, $x_n$, then $w_{\circ m+t} (x_i) \ne e$ for all $t\ge 1$.

\begin{cor} \label{cor:generaldynamics}
$k\ge 1$, take words $w_1(a_1, \dots, a_k)$, \dots $w_k(a_1, \dots, a_s)$ and suppose that for 
all $m$
the words $w_{j, \circ m} $ are not freely equivalent to identity, for all $j : 1 \le s$.
Then there exist a finite group $G$ such that for all $m\ge 1$ and all $j : 1 \le j \le  s$ the iterations $w^i_{\circ m} \ne e$ in $G$.

\end{cor}

{\it Proof of Theorem \ref{thm:borsapmol} and Corollary \ref{cor:generaldynamics}} 

We know that for all $m$, and hence in particular for $m=4s$ that the words $w_{j, \circ m} $ are not freely equivalent to identity, for all $j : 1 \le s$.

Consider $s$ matrices $M_j$ over  $\mathbb{Z}[x_{i_1,i_2, j}])$, where $i_1,i_2: 1 \le i_1, i_2 \le 2$, $j: 1 \le j \le s$ and $x_{i,j}$ are independant variables:
\[ M_j= \left( \begin{array}{ccc}
x_{1,1,j} & x_{1,2,j}  \\
x_{2,1,j} & x_{2,2,j}
 \end{array} \right) 
\]

Consider rational functions $R^{(n)}_{r_1,r_2,t}$ in $x_{i_1,i_2,j}$ , $i_1,i_2: 1 \le i_1, i_2 \le 4$, $j: 1 \le j \le s$ and $r_1: 1 \le r_1,r_2 \le 2$, $s: 1 \le s \le s$ which are the entries of $w_{\circ n, t}(m_1,\dots, m_s)$, $t: 1 \le t \le s$,
these rational functions are of the form 
$$
R^{(n)}_{r_1,r_2,t}= P^{(n)}_{r,t}/\prod _j (\det M_j)^{\alpha_j},
$$
 where $P^{(n)}_{r,t}$ are polynomials in $x_{i_1,i_2,j}$ with integer coefficients, and $\alpha_j$ are some integers.

We want to find a non-trivial solution over a finite field of the system of equations for some $n\ge 1$
\begin{equation}  \label{systemR}
R^{(n)}_{i_1,i_2,j} =x_{i_1, i_2,j}.
\end{equation}

To to this, we want to find a solution of the system of the equation
\begin{equation} \label{systemP}
P^{(n)}_{i_1,i_2,j}(x_{r,t}) =x_{i_1, i_2,j},
\end{equation}
$r \le 4$, $t \le s$,  where none  of the  two-times-two matrices $m_j$, $j: 1 \le j \le s$ is proportional to a diagonal matrix and satisfying $\det m_j = x_{1,j} x_{4,j} - x_{2,j} x_{3,j} \ne e$, for all $j: 1 \le j \le s$.

To do this it is sufficient to find, for some large power $Q$ of $q$ a  solution of the system of the equations, each $M_j$ has determinant not equal to zero, none of $M_j$ is  proportional to an identity  matrix, 
none of the coordinates of the $u$-th iteration ($u \le m$, $m$ is an appropriate function of $Q$, $n$, and $s$ )  of $w$ applied to $M_1, \dots, M_s$ is proportional to the identity matrix.
\begin{equation} \label{systemPQ}
P_{i_1, i_2,j} =x_{i_1, i_2,j}^Q 
\end{equation}
over a finite field of characteristics $q$.

\begin{rem} Let $W$ is a  word in $x_1$, \dots, $x_n$ which  is not freely equivalent to an empty word
Then the entry  $M_{1,2}$   in the upper-right corner of the matrix $M$ (which is a rational function in $x_1$, \dots, $x_n$) for the matrix
$$
M = W(m_1, \dots, m_2)
$$
is not equal to zero.
\end{rem}

\begin{rem} \label{re:nonzero}
Let $W_1$, $W_2$, \dots, $W_N$ are words in $x_1$, \dots, $x_n$ such that none of these words  is freely equivalent to an empty word.
Let $F_1$, \dots, $F_L$ are integer valued polynomials in $x_1$, \dots, $x_n$, each of $F_i$ is not identically zero.
Then there exists a finite field $\mathcal{K}$ and $x_1$, \dots, $x_n \in \mathcal{K}$ such that the  (upper-right) entry $M_{1,2,j}$ of the matrix
$$
M_j = W_j(M_1, \dots, M_2)
$$
is not equal to $0$, for each $j: 1\le j\le  N$ and such that
$F_j(x_1, x_2, \dots, x_n) \ne 0 $ for all $j \le L$.

\end{rem}

Now consider  $n=N=s$ , $W_{j}=w_{\circ 4s,j}$.
 From the assumption of the theorem we known that none of the words $W_j$ is freely equivalent to an empty word, and hence these words
satisfy the assumption of Remark \ref{re:nonzero}. Consider $M=s$ and $F_j = x_{1,j} x_{4,j} -x_{2,j} x_{3,j}$. From Remark  \ref{re:nonzero} we know that there exist a point $v_{i,j}$, $ i \le 4$, $j\le s$
in the image of the mapping corresponding
to $w_{\circ 4s,j}$ such that $ v_{1,j} v_{4,j} -v_{2,j} v_{3,j} \ne 0$ for all $j$ and such that $v_{2,j} \ne 0$ for all $j$.

Consider the polynomial in $x_{j,i}$, $j: 1 \le j \le s$, $i: 1 \le i \le 4$
$$
D =\prod_{j=1}^s x_{j,2}  \prod_{j=1}^s  (x_1 x_4 -x_2 x_3)
$$

The degree of $D$ is equal to $3s$, and $D(v_{i,j}) \ne 0$ for some $v_{i,j}$ in the image of the mapping corresponding to $w_{\circ 4s,j}$. Applying Theorem \ref{thm:borisovsapir}  for $n=4s$ and $d$ to be equal to the length
of $W$
we can conclude that there exists a solution $x_{i,j}$ for the system of the equations  \ref{systemR},  such that $D(x_{i,j}) \ne 0$
so far as $Q$ satisfies
$$
Q/3s > (4s)(4s+1) d^{16 s^2+1}
$$

\section{Open questions} \label{se:openquestions}

We recall again that our main interest are the words in the commutator, with total number of $x$ which is not zero (and with total number of $y$ which is not zero). 
Take  $w(x,y)$ is such that the total number $X$ of $x$ is not zero.
The total number of $x$ in the $n$-th iteration is equal $X^N$, and the total number in $w_{\circ n}(x)x^{-1}$ is $X^n-1$. So if $X\ne 2$, already without taking any iteration
$w(x,y) =x$ has a solution with $x\ne 0$ in an Abelian finite cyclic group, and if $X=2$, the equation for the second iteration $w(w (x,y),y)=x$ has an equation in a finite Abelian group
$\mathbb{Z}/2\mathbb{Z}$. For example, if we take $w(x,y)= yx^2y^{-1}$, then for the first iteration we obtain the solvable Baumslag Solitar group (so that the equation
does not have solution in finite groups with $x\ne e$), but for the second iteration we do obtain such solutions.

Given a word $w$, one can ask what is minimal $m$, which we denote by $m(w)$, such that $w_\circ(x)=x$ has a non-zero solution in a finite group. What is the minimal size $M(w)$ of a finite group which does not satisfy the iterated identity $w$.

Absence of (usual) identities in the class of finite quotients of a given group $G$ (for example absence of identities for all finite groups, or all finite nipotent groups etc for $G=F_m$) can be a corollary of residual finiteness
of $G$. One can make the statement quantative, by taking a word $w$ of length  $l$ and ask for 
 a minimal possible size of finite quotient of $G$ which does not satisfy $w$. Or a less stronger version: for a minimal possible size where $w(x_1, \dots, x_n) \ne e$ for a fixed finie set $x_1$, $x_2$, \dots, $x_n$ in $G$.
This notion, introduced by Bou-Rabee in \cite{bourabeedef} is  called {\it the normal residual finiteness growth function}, see also \cite{bourabeemcreynolds, bourabee7, buskin, kassabovmatucci},
and it is called {\it residual finiteness growth function} in  Thom \cite{thom} and Bradford  and Thom \cite{bradfordthom} (not to be confused with residual finiteness growth function in terminology
of \cite{bourabeeetal}, that measures the size of finite, not necessary normal subgroups, not containing a given element), who have proven the lower bound $\ge C n ^{3/2}/\log^{9/2+\epsilon}{n}$, which holds
for any $\epsilon>0$. Kassabov and  Matucci suggested in  \cite{kassabovmatucci} that the argument of Hadad \cite{hadad} can give a close upper bound for normal residual finiteness growth, function, namely $n^{3/2}$. 
A known upper bound so far is $n^3$ \cite{bourabeedef}, which is a corollary of the estimate for $SL(2, \mathbb{Z})$, using imbedding of a free group to this group.
The estimate of Bradford and Thom is a corollary of their result, stating that for any $\delta>0$ and all $n\ge 1$ there exists a word $w_n$ of length at most $n^{2/3} \ln^{3+\delta}(n)$ which is an identity in all
finite groups of cardinality at most $n$.
Now one can ask corresponding questions related to iterated identities. In particular, one can ask, what is the minimal length of a word $w_n$ which is an iterated identity in all finite groups of cardinality
at most $n$? Given a word $w$, we denote  by $NI(w)$ the minimal cardinality of a group $G$, such that $w$ is not an iterated identity in $G$ and by $PE(w)$ the the minimal cardinality of a group $G$ such that
$w_\circ{m}(g) =g$ has at least one non-identity solution in $G$, for some $m\ge 1$. It is clear that $PE(w) \ge NI(w)$. We also denote by $PE_d(n)$ and $NI_d(n)$ the maximum of $PE(w)$ and $NI(w)$, where
the maximum is taken over all words of length at most $n$ on $d$ letters, not freely reduced to an empty word. Finally, given a word $w$ we can ask what is the minimal $m$ such that 
$w_\circ{m}(g) =g$ has at least one non-idenity solution in some finite group?

Another question we can ask: what are possible classes of finite groups, where with the property that for any $w$, not freely equivalent to the identity,  there exists a group in this class which does not satisfy an iterated identity $w$. In particular, given a subset $\Omega \subset \mathbb{N}$, one can ask: for which subsets  $\Omega$, for any $w$, not freely equivalent to the identity,  there exists a group $G$, with the cardinality of $G$ belonging
to $\Omega$, such that $G$ does not satisfy the iterated identity $w$. We have seen in the proof of Theorem \ref{thm:noid} that it is sufficient to consider $SL(2,F_Q)$, for $Q$ which is a large power of a large prime $q$
and hence  the set $\Omega$ containing numbers $q^n-q$, for large enough $q$ and large enough $n$, has this property. We denote by $ \mathcal{O}_{int}$ the set of $\Omega \subset \mathbb{N}$ with the property above.
By $\mathcal{O}$ we denote the set of subsets $\Omega \subset \mathbb{N}$ such that for any word $w$, not freely equivalent to the identity, there exists a finite group $G$, of the cardinality belonging to $\Omega$ such that
$w$ is not an identity in $G$.
It is clear that  $\mathcal{O}_{int} \subset \mathcal{O}$  and that  $\mathcal{O}_{int} \ne\mathcal{O}$ since the set of powers of a given prime $p$ belongs to $\mathcal{O})$ for all $p$ and does not belong to
$\mathcal{O}_{int}$.

\end{document}